\documentclass[12pt]{amsart}
\usepackage{amssymb,cite}
\usepackage{amsmath,graphicx}
\begin{document}
\begin{center}
\includegraphics[scale=0.7]{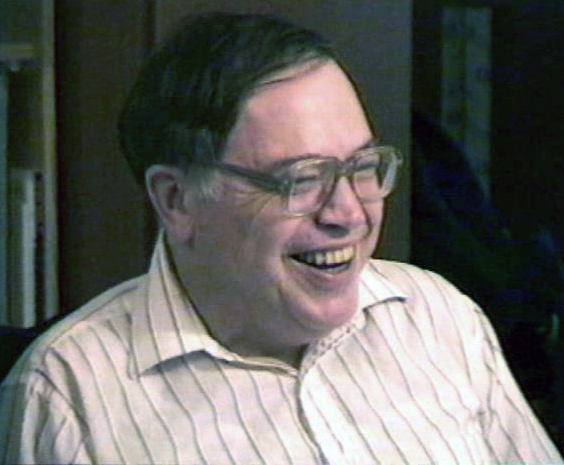}\\
Misha Shubin\\
1944 -- 2020
\end{center}

Mikhail Aleksandrovich Shubin, an outstanding mathematician, passed away after a long illness on May 13, 2020. Mikhail Aleksandrovich was born on December 19, 1944 in Kujbyshev (now Samara). He was brought up by his mother and grandmother. His mother, Maria Arkadievna, was an engineer at the State Bearing Factory, where she was hired in 1941 after having graduated from the Department of Mechanics and Mathematics of the Moscow State University (MSU). At that time, the factory was evacuated from Moscow to Kujbyshev. Maria Arkadievna had been working at the factory for many years as the Head of the Physics of Metals Laboratory. Later, she defended a PhD dissertation and moved to the Kujbyshev Polytechnical Institute, where she worked as an Associate Professor.

In his  school years, Shubin was mostly interested in music. He had an absolute pitch,  finished a music school, and was seriously thinking of entering a conservatory. However, while in high school, he developed an interest to mathematics, was successful in olympiads, and eventually decided to apply to the Department of Mechanics and Mathematics (mech-mat).

Mikhail Aleksandrovich was admitted to mech-mat of MSU in 1961. When the time came to choose an adviser, he became a pupil of Marko Iosifovich Vishik. After graduation, he was enrolled to the Graduate Program, and in 1969 defended his PhD Thesis. In the dissertation he derived formulas for computing the index of matrix-valued Wiener-Hopf operators. In particular, for the study of families of such operators, he had to generalize a theorem of Birkhoff stating that a continuous matrix-valued function
 $M(z)$ defined on the unit circle $|z|=1$ can be factored as
\begin{equation*}
M(z) =A_+(z) D(z) A_-(z),
\end{equation*}
where $A_+(z)$ and $A_-(z)$  are continuous and have analytic continuations to the interior of the unit circle and its exterior (infinity included), respectively; $D(z)$  is a diagonal matrix with the entries $z^{n_j}$  on the diagonal, with integer $n_j$. Shubin considered the problem about what happens when the matrix $M$ depends continuously on an additional parameter $t$. The Birkhoff factorization can not be made continuous in  $t$. Mikhail Aleksandrovich showed that if one relaxes the condition that $D(z)$ is diagonal, allowing it to be just triangular, then the factorization can be done continuously with respect to $t$.

Another paper \cite{70} from that period deserves to be mentioned. Let a family of closed subspaces of a Banach space depend analytically on a multi-dimensional complex parameter. This means that locally with respect to the parameter there exists a holomorphic projector-valued function, with the subspaces from the family being the ranges of the corresponding projections. One wonders whether there exists such a globally holomorphic projector-valued function. Shubin proved that if the domain of the parameter is  holomorphically convex (or a more general Stein complex space), the answer is positive. That required application of a recent result of L.~Bungart, which  generalized to an infinite-dimensional situation a difficult and not old at that time theorem by H.~Grauert claiming that holomorphic and topological classifications of vector bundles over Stein spaces coincide. This result of Shubin has become classical and extensively used in the operator theory. The paper was one of early manifestations of the broad range of Mikhail Aleksandrovich's knowledge.

The notion of a pseudodifferential operator was formalized in the mid-sixties, and the rapid development of the  micro-local analysis ensued. The micro-local analysis, one way or another became central to practically all works by Shubin. In collaboration with V.~N.~Tulovsky \cite{Tulov}, he studied the spectral asymptotics of an operator $P$ of order $m>0$ with a Weyl symbol in $R^n$. The spectrum of this operator is discrete, if the symbol satisfies the following conditions:
\begin{equation}\label{E:shubinclass1}
|D^\alpha p(y)  |\le C_\alpha (1+|y|)^{m-\rho | \alpha|},
\end{equation}
where $0<\rho\leq 1$,
\begin{equation}\label{E:shubinclass2}
|D^\alpha p(y)| \le  C_\alpha p(y)|y|^{-\rho|\alpha |},
\end{equation}
and
\begin{equation}\label{E:hypoell}
 p(y) \ge  a|y|^{m_0}
\end{equation}
 for large $|y|$, where the numbers  $a$  and $m_0$ are positive. Under some additional assumptions about the symbol, they proved that when $\lambda$ goes to infinity, the number of the eigenvalues that do not exceed $\lambda$ asymptotically equals   $(2\pi)^{-n}$ times the measure of the set $\{y \colon p(y)<\lambda \}$.   They also obtained a remainder estimate. The method of the proof (which has became popular since) was remarkable: they applied techniques of micro-local analysis to construct  approximate spectral projections, and then used variational methods.  The class of operators with symbols satisfying (\ref{E:shubinclass1}) is frequently called the Shubin class.

In the second half of seventies, Mikhail Aleksandrovich worked mostly on spectral theory of operators with almost periodic and random coefficients. He introduced a class of pseudodifferential operators with almost periodic symbols and studied spectral properties of selfadjoint operators from this class. He also proved existence of the density of states and of the Fermi energy for such operators. Let  $A$ be an elliptic selfadjoint almost periodic differential operator in $\mathbb{R}^n$  and  let $\Omega_k$ be a family of "good"  (e.g., with smooth boundary) increasing domains that exhaust the space and such that the ratio of their surface areas to the volumes goes to zero as $k$ increases (so-called F\"olner condition). Consider the operator $A$ in $\Omega_k$  with selfadjoint elliptic boundary conditions (e.g., Dirichlet).
Let   $E_j(k)$ be its eigenvalues taken in non-decreasing order and counted with their multiplicities. Let $N_k(E)$ be the number of these eigenvalues that are less than $E$. The density of states is the limit of $N_k(E)/|\Omega_k|$ when $k\to\infty$, and the Fermi energy $E^F(\rho)$ is the limit
\begin{equation*}
\frac{1}{p(k)} \bigl(E_1^{(k)} +\ldots + E_{p(k)}^{(k)}\bigr),
\end{equation*}
where  $k\to \infty$ and $p(k)/|\Omega_k|\to \rho$.

These works constituted his Doctor of Science (the degree in Russia roughly equivalent to Habilitation) Dissertation, which he defended in 1981 in Leningrad. There Shubin also developed functional calculus of pseudodifferential operators with almost periodic symbols, proved that the spectrum in different function spaces is the same, that the spectral $\zeta$ function admits analytic continuation, and derived a formula for the index \cite{PPOperators,CoincSpectra,PP_vonNeuman,PPF_PDEs,DensStates,PPIndex}. In particular, he used the $II_\infty$-factor introduced by Coburn, Moyer and Singer.

The methods used for studying operators with almost periodic coefficients are reminiscent of those that are known for operators with random coefficients. Following this connection, in collaboration with B.~V.~Fedosov \cite{Fed}, Shubin derived a formula for the index of elliptic operators in $\mathbb{R}^n$ with homogeneous random fields as coefficients.
Because these operators are not Fredholm in the standard sense, their index is defined in terms of traces in von Neumann algebras. Together with S.~M.~Kozlov \cite{Kozlov,Kozlov2}, they studied the problem of  equality of the spectra in different functional spaces.

Shubin has made also a significant contribution to the theory of operators with periodic coefficients. In the joint work with D.~Schenk \cite{Shenk,Shenk2} they derived a complete asymptotic expansion (when the energy goes to infinity) of the density of states for the Hill operator with a smooth potential. One of the difficulties here came from the fact that for a generic potential all gaps in the spectrum are open, and thus the expansion can not be differentiable. Only recently, the multi-dimensional analog of this result has been obtained by L.~Parnovski and R.~Shterenberg.

In 1982 E.~Witten introduced a deformation of the de Rham complex associated with a Morse function on a closed manifold, and used it for a new proof of the Morse inequalities. S.~P.~Novikov and M.~A.~Shubin \cite{Novikov} noticed that the Witten's deformation can be applied to the following situation.  Let $M$ be a closed manifold with an infinite fundamental group, and let  $\widetilde{M}$ be its universal covering space. Introduce a Riemannian metric on $M$  and lift it to $ \widetilde{M}$. Consider the Laplacian $\Delta_p$ acting on $p$-forms on $\widetilde{M}$, and let
$K_p(x,y)$ be the Schwartz kernel of the orthogonal projection onto the null-space of $\Delta_p$. The fundamental group  $\pi_1(M)$  acts on $\widetilde{M}$, and $K_p(x,y)$ is invariant with respect to this action. One can define $L^2$-Betti numbers as follows:
\begin{equation*}
\beta_p(M):=\int_F \operatorname{tr} K_p(x,x)\, dx,
\end{equation*}
where $F$ is a fundamental domain of the action of $\pi_1 (M)$ on $\widetilde{M}$.
Novikov and Shubin proved that if $f(x)$ is a Morse function on $M$, then the Morse inequalities remain valid if one replaces the usual Betti numbers by their $L^2$ analogs. They also introduced a new system of invariants of a manifold $М$ that was named later the Novikov--Shubin invariants. Let
$N_p(\lambda ,x,y)$ be the Schwartz kernel of the spectral projection for the operator  $\Delta_p$, and let
\begin{equation}
N_p(\lambda):=\int_F \operatorname{tr}N_p(\lambda, x,x)\, dx.
\end{equation}
The invariant $\alpha_p$ is the smallest number, for which
\begin{equation*}
N_p(\lambda )-\beta_p(M)=O(\lambda^{\alpha_p})
\end{equation*}
as $\lambda\to 0$. Later, Shubin and M.~Gromov showed \cite{Gromov} that $\alpha_p$ are homotopic invariants of $M$.  These papers with Novikov and with Gromov treat also a more general case of a representation of the fundamental group $\pi_1(M)$ in a $II_1$-factor.

The classical Riemann--Roch theorem can be interpreted as a formula for the change of the index of the Cauchy--Riemann operator that occurs when one enforces a divisor of zeros and singularities.  Some generalizations of this theorem to elliptic operators were obtained in the works by V.~G.~ Maz'ya and B.~A.~Plamenevski and N.~S.~Nadirashvili. M.~A.~Shubin, in collaboration with M.~Gromov, obtained far reaching generalization of these theorems for elliptic operators acting on sections of vector bundles over closed manifolds, as well as for the ``divisors'' that are not necessarily discrete. We will formulate one of their results. Let $A$ be an elliptic differential operator acting on sections of a vector bundle  $E$ over a closed manifold  $M$;  $E^*$ is the dual bundle to $E$ and $A^*$ is the transposed operator to $A$. A rigged divisor is a quadruple $(D^+,L^+;D^-,L^-)$ where $D^+$, $D^-$ are disjoint nowhere dense closed sets in $M$, and  $L^+$, $L^-$ are finite-dimensional spaces of distributions with the values in $E$ and $E^*$ and with supports lying in  $D^+$  and $D^-$ respectively. By $L(\mu, A)$  we denote the space of smooth sections of $E$ over $M \backslash D^+$ that can be extended to distributions $u$ over $M$ with values in $E$ such that $u$ is orthogonal to $L^-$, and  $Au\in L^+$ . Then
\begin{equation}
\operatorname{dim} L(\mu, A) =\operatorname{ind} A + \operatorname{deg}_A(\mu) + \operatorname{dim} L(\mu^{-1}, A^*).
\end{equation}
Here  $\mu^{-1}:=  (D^-,L^-; D^+,L^+)$ and  $\operatorname{deg}_A(\mu)$ is some ``degree'' of the rigged divisor (usually explicitly computable), which we will not define here. Note that if one writes the above equality as
\begin{equation*}
\operatorname{dim} L(\mu, A) - \operatorname{dim} L(\mu^{-1}, A^*) =\operatorname{ind}  A + \operatorname{deg}_A \mu, A,
\end{equation*}
it can be interpreted as the perturbation of the index occurring when one enforces a divisor of ``zeros and singularities.''

Later, Shubin obtained an  $L^2$-version of this result for non-compact regular coverings of compact manifolds.

In the second half of the 90s and in the beginning of 2000s, Shubin extensively studied, both on his own and in collaboration with V.~G.~Maz'ya and V.~A.~Kondratiev, criteria of discreteness of the spectrum and of essential self-adjointness for Schr\"odinger operators with electric and magnetic potentials \cite{KMS_magnetic,KMSgauge,Noncomp3,Kondr2,Maz'yaDiscr}. These papers are related to the following result of A.~M.~Molchanov that was obtained in 1953. Let $H=-\Delta +V(x)$ be a Schr\"odinger operator in  $\mathbb{R}^n$  with a locally integrable potential $V(x)$. A.~M.~Molchanov showed in 1953 that the following condition is equivalent to the discreteness of the spectrum of the operator $H$: for every $d>0$
\begin{equation*}
\inf \int _{Q_d \setminus F} V(x) \,dx  \to \infty
\end{equation*}
when $Q_d\to \infty$ .
Here $Q_d$ is a cube of side $d$, the sides of which are parallel to the co-ordinate axes, and the infimum is taken over compact sets $F$ such that their Wiener capacity  $\operatorname{cap} F$  does not exceed  $\gamma\operatorname{cap}\overline{Q_d}$ (such sets are called $\gamma$-negligible), and $\gamma$ is small enough. Giving an answer to I.~M.~Gel'fand's question, Maz'ya and Shubin showed that every $\gamma$ such that $0<\gamma <1$ works. Moreover, one can make $\gamma$ depend on $d$ in such a way that
\begin{equation*}
\limsup_{d\to 0} \gamma(d)d^{-2}=\infty.
\end{equation*}

For operators with a magnetic potential, the discreteness of the spectrum criteria were obtained in terms of an ``effective'' scalar potential introduced for that purpose.

In the paper  ``Can one see the fundamental frequency of a drum?'' \cite{drum} Maz'ya and Shubin obtained a beautiful and far reaching generalization of the following result by W.~K.~Heiman: let $\Omega$ be a simply connected domain in $\mathbb{R}^2$ and let  $r(\Omega)$ be the maximal radius of a circle inscribed in $\Omega$ (``inner radius''). Then the smallest eigenvalue of the Dirichlet Laplacian in  $\Omega$  satisfies the estimate
\begin{equation*}
cr^{-2}(\Omega) \le \lambda_1(\omega)\le Cr^{-2}(\Omega)
\end{equation*}
with absolute constants  $c$ and $C$. This result does not hold for non-simply connected domains (drilling a large number of small holes in a domain gives a counter-example), and it also does not hold in higher dimensions. Maz'ya and Shubin showed that if one changes the definition of the inner radius appropriately, then one can extend the result to arbitrary domains in arbitrary dimension $n$. Namely, let $0<\gamma<1$.
Let $r_\gamma(\Omega)$ be the maximal radius of a closed ball $B_r$ such that $B_r\setminus \Omega$  is $\gamma$-negligible in $B_r$. Then
\begin{equation*}
c(n,\gamma) r_\gamma^{-2}(\Omega) \le \lambda_1(\omega)\le C(n,\gamma)r_\gamma^{-2}(\Omega).
\end{equation*}

Mikhail Aleksandrovich had also worked in a variety of other areas of mathematics. One should mention in particular his undeservedly little known paper on pseudo-difference operators and on the estimates of discrete Green's functions \cite{Difference} (see also \cite{DiscrMagn,Sunada}). Nowadays, discrete spectral problems attract a lot of attention, and Shubin's technique should be extremely useful. It is impossible to adequately describe in this article his numerous papers (both single authored and joint  with his students and collaborators) devoted to the spectral theory of operators on non-compact manifolds and on Lie Groups (in particular, see  \cite{TransvEllipt,Brav,Noncomp,Noncomp2}).  He had also worked on completely integrable equations \cite{Kapp1,Kapp2}, applications of non-standard analysis \cite{Nonstand}, and other problems. The probably incomplete list of his publications in MathSciNet contains 135 items.

Shubin has written several remarkable books. Published in 1978, his book "Pseudodifferential operators and spectral theory" \cite{PseudoBook} instantly became one of the main textbooks on microlocal analysis. Now, after more than forty years, it is still very popular. His joint monograph with F.~A.~Berezin "The Schr\"odinger Equation" \cite{Berezin} has been very influential. In that book, which was finished after Felix Alaksandrovich Berezin passed away, many topics, e.g. the path integrals, are treated in a rather unique way. In 2001 Mikhail Aleksandrovich wrote the textbook ``Lectures on the equations of mathematical physics'' \cite{UrMatFiz}, which was based upon the lectures that he gave at the Moscow State University. Somewhat extended English translation of this book, edited by M.~Braverman, R.McOwen, and P.~Topalov, was recently published by the AMS  \cite{PDEBook}. In addition to numerous expository articles that he wrote for Russian Math. Uspehi, one should also mention extensive surveys that he wrote in collaboration with other authors for the VINITI publications  (for example,  \cite{Egorov,Roz}). He also edited Russian translations  of the fundamental books by F.~Treves and L.~H\"ormander.

Mikhail Aleksandrovich started participating in the famous seminar of I.~M.~Gel'fand in October of 1964, when he was a junior undergraduate student. He had been taking notes, practically with no interruptions, for 25 years, and preserved all of them. Now, with the financial support from the Clay Institute, these notes are available on the internet as the  ``Gelfand Seminar notes''  file \cite{GelfSeminar}.
Gelfand's seminar was the center of mathematical life in the USSR, and Shubin's notes have become a  unique testament  of that years.

Mikhail Aleksandrovich was a remarkable lecturer and teacher. His lectures were extremely well prepared, with all small wrinkles ironed out. In a rather short period if time, he would lead students to non-trivial problems. He had advised on about 20 PhD dissertations. In addition to his PhD students, he strongly influenced many other young mathematicians. This was not limited to undergraduate and graduate students. He had worked a lot with high school students in summer math camps; as a result he published some of his lecture notes there and a small book \cite{MatPros} in the ``Mathematics Education'' series. Starting 1969, Shubin worked at the Division of Differential Equations of the Department of Mechanics and Mathematics of the Moscow State University. In 1992 he became a Professor, and then Distinguished Professor of the Department of Mathematics of the Northeastern University in Boston, USA.

Mikhail Aleksandrovich was a remarkable, broadly educated, cheerful, and friendly person. In the seventies and eighties he was helping several Soviet mathematicians who had difficulties due to political reasons. He also fought against unfair admission practices that were taking place at the Department of Mathematics of the MSU.

The memory of Mikhail Aleksandrovich Shubin, a mathematician, teacher, and a human being, will remain forever in the hearts of his numerous colleagues, students, and friends.

\begin{flushright}
\emph{M.~Braverman, B.~M.~Buchshtaber,  M.~Gromov, V.~Ivrii, Yu.~A.~Kordyukov, P.~Kuchment, V.~Maz'ya, S.~P.~Novikov, T.~Sunada,
L.~Friedlander, A.~G.~Khovanskii}
\end{flushright}

\end{document}